\hoffset -4pt
\baselineskip 13pt
\magnification 1850

\input amssym.def
\input amssym.tex


\font\bmat=cmbxti10 at 14pt
\font\csc=cmcsc10
\font\title=cmbx12 at 14pt

\font\teneusm=eusm10    
\font\seveneusm=eusm7  
\font\fiveeusm=eusm5    
\newfam\eusmfam
\def\eusm{\fam\eusmfam\teneusm}
\textfont\eusmfam=\teneusm 
\scriptfont\eusmfam=\seveneusm
\scriptscriptfont\eusmfam=\fiveeusm


\def\Re{{\rm Re}\,}

\def\sgn{{\rm sgn}\,}
\def\txt#1{{\textstyle{#1}}}

\def\r#1{{\rm #1}}
\def\e#1{{\eusm #1}}
\def\varGamma{{\mit\Gamma}}
\def\B#1{{\Bbb #1}}

\def\sgn{{\rm sgn}}
\centerline{\title Spectral Analysis}
\medskip
\centerline{\title of the Zeta and {\bmat L}$\,$-Functions}
\bigskip
\centerline{\csc Yoichi Motohashi}
\footnote{}{\hskip-0.7cm Talk at the workshop
 `{\it L-functions of automorphic forms and related problems$\,$}'
March 2012, Tokyo University. }
\medskip
\noindent
{\bf 1. Aim of talk}
\par
\noindent
We will show an instance of applications of the theory of
automorphic representations to a genuinely traditional problem 
in the theory of the zeta and allied functions. We shall 
restrict ourselves
to very basic issues and results, 
because of the purpose of the workshop.
\medskip
\noindent
{\bf 2. The original problem}
\par
\noindent
Bounding the zeta-function $\zeta(s)$, $s=\sigma+it$, especially
on the critical line $\sigma={1\over2}$, 
is one of the central problems
in analytic number theory, as it is indispensable in the theory of the
distribution of prime numbers and more centrally in investigating
the distribution of complex zeros of the zeta-function. 
This is the same as treating 
non-trivially the zeta-sum
$$
\sum_{N<n\le 2N}n^{it},
$$
where $N, t>0$ are large and independent of each other. In other
words, a massive cancellation among the 
waves $\{n^{it}\}$ is expected and to be detected.
\vfill
\break
\noindent
{\bf 3. Divide and conquer}\footnote\dag
{Appellation by M.N. Huxley.}
\par
\noindent
The first significant contribution was achieved by H. Weyl (1921).
It was greatly
improved by J.G. van der Corput (1921--37) using a tool
from harmonic analysis (Poisson's sum formula).
\par
Simplifying the story, we
chop up the zeta-sum into pieces, and apply 
Cauchy's inequality, getting the expression
$$
\sum_{n\in\B{N}}
w(n)\Bigg|\sum_{0<m\le M}(m+n)^{it}\Bigg|^2,
$$
with an appropriate test function $w$. We then
expand the squares, take the sum over $n$ inside,
and appeal to Poisson's sum formula in handling
{\it non-diagonal\/} parts.\footnote\ddag{Thus follows 
the famous subconvexity
bound of van der Corput: $\zeta({1\over2}+it)\ll t^{1/6}\log t$, 
$t\ge2$.} 
\medskip
The resulting non-trivial bounds were exploited
in the theory of the zero-free region and 
the zero-density concerning the
zeta-function, and yielded
eventually G. Hoheisel's 
detection (1930) of prime numbers in {\it short\/}
intervals, without any hypothesis like Riemann's. It was the dawn of
the modern theory of the zeta-function and the distribution of
prime numbers.
\vfill
\break
\noindent
{\bf 4. The most fundamental problem in ANT}
\par
\noindent
This is to discuss the sum
$$
\sum_n f(n).
$$
\medskip
\noindent
{\bf 5. Lifting to higher dimension}
\par
\noindent
This is essentially due to F.V. Atkinson (1949) and can be
regarded as a generalisation of Weyl's idea. Thus, we
square the last sum and consider instead the two dimensional
sum
$$
\sum_m\sum_nf(m,n),
$$
where $f$ has a new specification (the same convention applies
to later steps). Then, we classify the summands 
according to the deviation
from the symmetric axis $m=n$. The sum
$$
\sum_m\sum_n f(n, n+m)
$$
emerges.
\bigskip
\noindent
{\bf 6. A surprise}
\par
\noindent
Then, we apply Poisson's sum formula to the last inner sum,
assuming $f$ is appropriate to do so. Often a better result
comes out. Well, it is better than what can be
achieved by applying Poisson's sum formula
to the original double sum. This mysterious effect of
{\it divide and conquer\/} was 
observed by van der Corput in an explicit way.
\vfill
\break
But, why is such a mystery possible? 
We are not claiming any full excavation of this under-sea 
{\it hoard\/} but only try to reveal
a group structure and indicate that 
an ocean of orthogonal waves embraces
the mystery. 
\medskip
\noindent
{\bf 7. Seeking axis of symmetry}
\par
\noindent
It is reasonable to surmise that the axis $m=n$ in the original
double sum is {\it not\/} in general the central line of symmetry. 
Also, we note that the application of Poisson's 
sum formula is the same as
looking for hidden orthogonalities. As much as, but not too much,
orthogonalities are to be exploited. The situation reminds us 
of the circle method of Hardy--Littlewood--Ramanujan: To seek
an {\it optimal\/} dissection.
\medskip
\noindent
{\bf 8. Further lifting}
\par
\noindent
Being in the NT community, we try the line $km=ln$ with an arbitrary 
pair of non-zero integers $k,l$; that is, turning the line $m=n$ to
$m=(l/k)n$ in order to achieve the best focusing.
The Farey sequence is visible here, in an analogous 
context as in the circle method.
Since we do
adjusting inside a view finder, weights are to be attached to
each pair $k,l$ in order to
see a particular set of points $\{l/k\}$ better than the rest
so that the control is kept in our hands.
\vfill
\break 
We are thus lead to a division of general quadruple sums:
$$
\eqalign{
&\sum_{k,l,m,n}f(k,l,m,n)\cr
=&\Bigg\{\sum_{km=ln}+\sum_{km>ln}+\sum_{km<ln}\Bigg\}
f(k,l,m,n)
}
$$
\medskip
\noindent
{\bf 9. Rendering with matrices}
\par
\noindent
It occurred to us, well more than 2 decades ago\footnote\dag
{Proc.\ Amalfi Conference 1989, ed.\ E. Bombieri et al,
Salerno Univ., 1992,
pp.\ 325--344.},
that the last identity could be better expressed
in terms of matrices: 
\medskip
With
$M=\left({k\atop n}{l\atop m}\right)\in\r{M}_2(\B{Z})$,
$$
\eqalign{
&\sum_{k,l,m,n}f(k,l,m,n)=\sum_{M}f(M)\cr
=&\Bigg\{\sum_{\det M=0}
+\sum_{\det M>0}+\sum_{\det M<0}\Bigg\}
f(k,l,m,n).
}
$$
\bigskip
\noindent
{\bf 10. Hecke and $\rm{SL}_2$}
\par
\noindent
We call the first sum on the right {\it the Ramanujan term\/}
with a good reason. Perhaps, {\it Rankin--Selberg\/} is to be
attached as well. The second and the third sums can be
associated with E. Hecke; perhaps, with {\it
Hurwitz--Mordell\/} also.
\vfill
\break
Thus we have
$$
\eqalign{
&\sum_{\det M>0}f(M)=\sum_{n>0}
(T(n)F)(1),\cr
F(\r{g})&=\sum_{M\in\r{SL}(2,\B{Z})}f(M\r{g}),\quad
\r{g}\in\r{SL}(2,\B{R}),
}
$$
where $\{T(n)\}$ are Hecke operators associated with the full
modular group.
\medskip
\noindent
{\bf 11. Spectral decomposition}
\par
\noindent
Assuming $f$ be sufficiently smooth and of rapid decay, we may
apply, to the {\it Poincar\'e series\/} $F(\r{g})$, the spectral
decomposition of the Hilbert space 
$L^2(\varGamma\backslash\r{G})$ or rather that
of automorphic representations occurring
there. Here $\r{G}=\r{SL}(2,\B{R})$ and $\varGamma=
\r{SL}(2,\B{Z})$. But a more convenient choice is
$$
\r{G}=\r{PSL}(2,\B{R}),\quad \varGamma=
\r{PSL}(2,\B{Z}),
$$ 
as we may assume that $f$ is even.
We shall keep this convention hereafter.
\bigskip
\noindent
{\bf 12. Upshot, in fact an intermezzo}
\par
\noindent
Applying the above consideration to the integral transform
$$
{\cal Z}_2(\zeta^2, g)=\int_{-\infty}^\infty\big|\zeta
\big(\txt{1\over2}+it\big)\big|^4g(t)dt,
$$
of a nice weight $g$, we are able to achieve a complete and
explicit spectral decomposition in terms of the
spectral resolution of
the Casimir operator acting over 
$L^2(\varGamma\backslash\r{G})$:
$$
\eqalign{
{\cal Z}_2(\zeta^2,g)&={\cal M}(g)+
\sum_V|\varrho_V(1)|^2H_V^3\big(\txt{1\over2}\big)
\Theta(\nu_V,g)\cr
&+{1\over4\pi i}\int_{(0)}
{|\zeta({1\over2}+\nu)|^6\over
|\zeta(1+\nu)|^2}\Theta(\nu,g)d\nu.
}
$$
There will be no need to explain the notation in detail.
We shall be a little more precise later. 
\medskip
This was achieved in a series of our works 
(1989--1997)\footnote\dag{Cambridge Tracts in Math., 
{\bf127}.}.
Initially we used the spectral theory of sums of
Kloosterman sums due to R.W. Bruggeman (1978) 
and N.V. Kuznetsov (1977--81). However, later their results were
dispensed with in our joint work\footnote\ddag{Crelle {\bf 579}
(2005), 75--114}  with  Bruggeman.
\par
The formula gave rise to a variety of new
facts on the zeta-function which had been unattainable before.
Main applications were done in our joint works 
with A. Ivi\'c concerning the plain fourth power mean
$$
\int_{-T}^T\left|\zeta\big(\txt{1\over2}+it\big)\right|^4dt,
$$
to which a part of Ivi\'c's talk at this workshop is devoted.
\vfill
\break
\noindent
{\bf 13. Real issues}
\par
\noindent
The spectral decomposition for ${\cal Z}_2(\zeta^2,g)$
is only the top tip of a great iceberg. It is beautiful
but is just the very beginning of a  
story to be unfolded in the future.
\par
There are at least two obvious directions 
in which we should go deeper:
$$
{\cal Z}_k(\zeta^2, g)=\int_{-\infty}^\infty\big|\zeta
\big(\txt{1\over2}+it\big)\big|^{2k}g(t)dt,
$$
and
$$
{\cal Z}_k(L_V, g)=\int_{-\infty}^\infty\big|L_V
\big(\txt{1\over2}+it\big)\big|^{k}g(t)dt,
$$
where, as is to be made precise later,
$L_V(s)$ is the $L$-function associated with an irreducible
cuspidal automorphic representation occurring in 
$L^2(\varGamma\backslash\r{G})$.
\par
We shall discuss only ${\cal Z}_2(L_V, g)$, 
since both for any integral $k\ge3$
belong presently to a terra incognita;
thus here is a challenging problem, whose resolution will
yield fundamental changes throughout ANT.
On the other hand, ${\cal Z}_2(L_V, g)$
is a natural extension of
${\cal Z}_2(\zeta^2,g)$. `Natural', because $\zeta^2$ 
corresponds to the Eisenstein series associated with $\varGamma$,
which is a kind of automorphic form. However, it is by no
means a ready-made extension, since the above scheme developed
for ${\cal Z}_2(\zeta^2,g)$ does not yield any incision,
especially when
dealing with $V$ in the unitary principal series, i.e., the irreducible
subspace generated by repeated applications of
Maass derivatives to a particular real-analytic
cusp-form on the hyperbolic upper half-plane.
\vfill
\break
\noindent
{\bf 14. An arithmetic--analytic issue}
\par
\noindent
In the Weyl--van der Corput setting, the function $f$
needs to be `smooth' actually. But in general $f$ is an
arithmetic function, and what is essential is to have a
sufficiently detailed asymptotical result on
the {\it shifted convolution\/}
$$
\sum_{n\in\B{N}} f_1(n)f_2(n+m)W(n/m),\quad m>0,
$$
with a test function $W$. Thus:
$$
\hbox{\it The Poisson sum formula is to be replaced }
\atop\hbox{\it by something more
arithmetic.}
$$
\medskip
When both $f_j$ are divisor functions, this is called an
{\it additive divisor problem/sum\/}
$$
\sum_{n\in \B{N}} d(n)d(n+m)W(n/m),
$$
to which corresponds
the mean value ${\cal Z}_2(\zeta^2,g)$.
This sum looks undisguised, but it is in fact a deep problem which has
an essential relation with the Ramanujan conjecture on the
size of Hecke eigenvalues. We need to appeal either to the
spectral theory of sums of Kloosterman sums or to
a careful construction of Poincar\'e series. Here the key-point
is that the divisor function has an inner-structure $\sum_{a|n}1$
which can be effectively and readily 
exploited.\footnote\dag{Our result on this is one of the 
essential implements
to prove a uniform subconvexity
bound for $L_V(s)$, which is a wide extension of
van der Corput's exponent $1\over6$ for $\zeta$.
See M. Jutila and Y.M: Acta Math., {\bf195} (2005), 61--115.
}
\medskip
Then, what will happen, if $f_j(n)$ are Fourier
coefficients of cusp forms? There does not seem to
exist any corresponding inner-structure. This problem was
posed by A. Selberg (1965), and was only recently
resolved by ourselves\footnote\ddag{Proc.\ Japan Acad., 
80{\bf A} (2004), 28--33.}
 appealing to the Kirillov model of 
automorphic representations
occurring in $L^2(\varGamma\backslash\r{G})$. 
We are going to indicate salient points of our idea.
\medskip
\noindent
{\bf 14. Normalisation}
\par
\noindent
We have the Iwasawa co-ordinate system
$$
\r{G}=\r{N}\r{A}\r{K}\ni \r{g}=\r{n}[x]\r{a}[y]\r{k}[\theta],
$$
with $\r{n}[x]=\left[{1\atop}{x\atop1}\right]$,
$\r{a}[y]=\left[{\sqrt{y}\atop}{\atop1/\sqrt{y}}\right]$,
$\r{k}[\theta]=\left[{\hfill\cos\theta\atop-\sin\theta}
{\sin\theta\atop\cos\theta}\right]$. The Casimir operator
is
$$
\Omega=-y^2(\partial_x^2+\partial_y^2)
+y\partial_x\partial_\theta.
$$
The (orthonormal) spectral decomposition is rendered as
$$
\displaystyle{L^2(\varGamma\backslash\r{G})
={\Bbb C}\cdot1\oplus
{}^0\!L^2(\varGamma\backslash\r{G})
\oplus{}^e\!L^2(\varGamma\backslash\r{G})}\atop
\displaystyle{{}^0\!L^2(\varGamma\backslash\r{G})
=\oplus V,\quad V=\mathop\oplus_\ell 
{\Bbb C}\lambda_V^{(\ell)}.}
$$
The subspace ${}^e\!L^2(\varGamma\backslash\r{G})$ 
is generated by
integrals of Eisenstein series, the details of which is omitted.
The set $\{\lambda_V^{(\ell)}\}$ is a complete 
orthonormal system of the cuspidal subspace
${}^0\!L^2(\varGamma\backslash\r{G})$, whose irreducible
subspaces are $\{V\}$; and
$$
\Omega\lambda_V^{(\ell)}=\big(\txt{1\over4}+\nu_V^2\big)
\lambda_V^{(\ell)},\quad \partial_\theta\lambda_V^{(\ell)}
=2i\ell\lambda_V^{(\ell)}.
$$
Also we have the Fourier expansion
$$
\eqalign{
\lambda_V^{(\ell)}(\r{g})=
\bigg|\pi^{-2\nu_V}{\Gamma\big(|\ell|+\nu_V
+{1\over2}\big)
\over\Gamma\big(|\ell|-\nu_V+{1\over2}\big)}
\bigg|^{1/2}\cr
\times\sum_{n\ne0}{\varrho_V(n)\over\sqrt{|n|}}
\e{A}^{\sgn(n)}\!\phi_\ell(\r{a}[|n|]\r{g},\nu_V).
}
$$
Here $\phi_\ell(\r{g},\nu)=y^{1/2+\nu}\exp(2i\ell\theta)$,
and the Jacquet transform
$$
\displaystyle{\e{A}^\delta\!\phi(\r{g})=\int_{-\infty}^\infty 
e(-\delta\xi)\phi(\r{w}\r{n}[\xi]\r{g})d\xi,}\atop{ 
e(\xi)=\exp(2\pi i\xi),
\quad \delta=\pm,\quad \r{w}=\r{k}
\big[\txt{1\over2}\pi\big].}
$$
This normalisation allows us to regard the Fourier coefficients
of cusp forms as being inherent to each $V$ but not the respective
forms.\footnote
\dag{See our exposition: arXiv:1112.4226v1[math.NT]}
\par
With this, we define
$$
L_V(s)=\sum_{n=1}^\infty\varrho_V(n)n^{-s},
\quad\Re s>1.
$$
The integral ${\cal Z}_2(L_V,g)$ is expressible in terms
of a generalised additive divisor problem/sum or 
the shifted convolution sum
$$
\sum_{n=1}^\infty
\varrho_V(n)\overline{\varrho_V(n+m)}W(n/m).
$$
The main problem is how to generate this sum by means of
the automorphic forms $\{\lambda_V^{(\ell)}\}$. 
\bigskip
\noindent
{\bf 15. Discrete series} 
\par
\noindent
If $V$ belongs to the discrete series, then $\nu_V=k-{1\over2}$
with a $k\in\B{N}$, and $\e{A}^\delta
\phi_\ell(\r{g},\nu_V)$ is essentially the exponential function.
Since the exponential function is a kind of additive character,
and we are done: In fact we may use the trivial but very basic
identity
$$
\int_0^1\exp(2\pi inx)\overline{\exp(2\pi i(n+m)x)}dx
=\delta_{m,0}.
$$
With this we may pick up $\{\varrho_V(n)
\overline{\varrho_V(n+m)}\}$ by considering the Fourier
expansion for $|\lambda_V^{(0)}|^2$.
\medskip
\noindent
{\bf 16. Unitary principal series}
\par
\noindent
If $V$ belongs to the unitary principal series, then $\nu_V\in i\B{R}$,
and $\e{A}^\delta\phi_\ell(\r{g},\nu_V)$ is essentially a
Whittaker function which is a generalisation of the $K$-Bessel
function (the Kelvin function)
and does not admit any property like an additive
character.
\medskip
{\it Then, is there any element in $V$ whose generic 
Fourier coefficient
is a product of $\varrho_V(n)$ and a function that admits an
additive property {\it similar\/} to the exponential
function?\/}
\medskip
\noindent
That is, we ask naively:
{\it  Is there an element in $V$ whose outward
appearance comes very close
to the discrete series?\/}
\bigskip
\noindent
{\bf 17. Kirillov map/model}
\par
\noindent
Our idea is to use the Kirillov model to find such an
automorphic form inside $V$. By the Kirillov map 
we mean the correspondence\footnote\dag
{loc.\ cit.}:
$$
\e{K}:\lambda^{(\ell)}_V\mapsto
\e{A}^{\sgn(u)}\phi_\ell(\r{a}[|u|],\nu_V),\quad
u\in \B{R}^\times.
$$
This extends linearly to the whole of $V$ and
becomes a unitary and surjective map between 
$V$ and $L^2(\B{R}^\times, d^\times\!/\pi)$,
with $d^\times u=du/|u|$. 
The right action $r$ of $\r{G}$ 
over the space $V$ is realised faithfully 
in $L^2(\B{R}^\times, d^\times\!/\pi)$:
$$
r_\r{h}\lambda_V^{(\ell)}\mapsto\e{A}^{\sgn(u)}(r_\r{h}
\phi_\ell)(\r{a}[|u|],\nu_V),\quad \r{h}\in\r{G}.
$$
\medskip
\noindent
{\bf 18. Inverse map}
\par
\noindent
With this, we pick up a function $\omega\in L^2(\B{R}^\times,
d^\times\!/\pi)$ which 
admits the necessary additive
property, and exploit $\e{K}^{-1}\omega$. For instance, we put
$$
\omega(u)=\cases{u^{\alpha+1/2}\exp(-2\pi u) & $u>0$,\cr
\hfil 0& $u\le0$,}
$$
with a sufficiently large $\alpha$ to gain a rapid decay on
the boundary. Then the function $(\e{K}^{-1}\omega)(\r{g})$ is
in the space $V$, and we have the expansion
$$
\eqalign{
&(\e{K}^{-1}\omega)(\r{n}[x]\r{a}[y])\cr
&=y^{\alpha+1/2}
\sum_{n>0}\varrho_V(n)n^\alpha\exp(2\pi in(x+iy)).
}
$$
This is similar to holomorphic cusp forms (i.e., elements in the 
discrete series). Then we consider the function
$$
|(\e{K}^{-1}\omega)(\r{g})|^2\in L^2(\varGamma\backslash
\r{G}),
$$
whose Fourier expansion involves $\{\varrho_V(n)
\overline{\varrho_V(n+m)}\}$ in a nice way.
The spectral decomposition of 
$|(\e{K}^{-1}h)(\r{g})|^2$ can be done explicitly
by using again the Kirillov map. We are essentially done.
{\it As a matter of fact\/}, we need to be a little more careful
in choosing the seed function $\omega$, 
but this point appears to be
immaterial in our present discussion.
\medskip
We have thus established the complete spectral
expansion\footnote\dag{Proc. Japan Acad., 83{\bf A} 
(2007), 73--78.} 
$$
\eqalignno{
{\cal Z}_2(L_V,g)&={\cal M}(g;V)
+\Re\Bigg\{\!\sum_U \varrho_U(1)
H_U\!\left(\txt{1\over2}\right)\Theta(\nu_U, g;V)\cr
+&{1\over4\pi i}\int_{(0)}\!{|\zeta\left({1\over2}+\nu\right)|^2
\over|\zeta(1+2\nu)|^2}L_{V\otimes V}
\!\left(\txt{1\over2}+\nu\right)\!\Theta(\nu, g;V)
d\nu\!\Bigg\}.
}
$$
This is an exact extension of our spectral expansion
for the fourth moment of the Riemann zeta-function. Here
$V$ is a particular irreducible representation while
$U$ runs over all cuspidal irreducible representations.
It should be stressed again that this result is {\it beyond\/}
the reach of the spectral theory of sums of Kloosterman sums. 
It belongs genuinely to
the theory of automorphic representations, that is, it appears
to us that without representation theory
of Lie groups the result would be very hard to achieve, if not
impossible. 
\medskip
As an application, we are currently developing a 
detailed {\it quantitative\/} analysis
of the plain mean value
$$
\int_{-T}^T\big|L_V\big(\txt{1\over2}+it\big)\big|^2dt.
$$
We have already obtained assertions analogous to those
on the fourth moment of the zeta-function. Our effort is
now focused to their uniformalisation with respect to $V$.
\medskip
Our way of using the Kirillov model has been 
exploited recently by V. Blomer and
G. Harcos\footnote\dag{arXiv: 0703246v1[math.NT]} 
and by others in investigating the problem of 
{\it uniformly\/} bounding various shifted convolution sums.
Also it is worth noting
that our idea can be extended to higher dimensional
situation in a straightforward way, though the analysis becomes
extremely complicated.

\bye